\documentclass[11pt]{amsart}
\usepackage{amsmath,amssymb,latexsym,esint,cite,mathrsfs}
\usepackage{verbatim,wasysym,cite}
\usepackage[left=2.8cm,right=2.8cm,top=2.9cm,bottom=2.9cm]{geometry}
\usepackage[latin1]{inputenc}
\usepackage{microtype}
\usepackage{color,enumitem,graphicx}
\usepackage[colorlinks=true,urlcolor=blue, citecolor=red,linkcolor=blue,
linktocpage,pdfpagelabels, bookmarksnumbered,bookmarksopen]{hyperref}
\usepackage[hyperpageref]{backref}
\usepackage[english]{babel}

\newtheorem{theorem}{Theorem}[section]
\newtheorem{proposition}[theorem]{Proposition}
\newtheorem{corollary}[theorem]{Corollary}
\newtheorem{lemma}[theorem]{Lemma}
\newtheorem{remark}[theorem]{Remark}
\newtheorem{definition}[theorem]{Definition}

\newcommand{\bcl}{\begin{center}}
\newcommand{\ecl}{\end{center}}
\newcommand{\brl}{\begin{right}}
\newcommand{\erl}{\end{right}}
\newcommand{\ben}{\begin{enumerate}}
\newcommand{\een}{\end{enumerate}}
\newcommand{\barr}{\begin{array}}
\newcommand{\earr}{\end{array}}
\newcommand{\btab}{\begin{tabular}}
\newcommand{\etab}{\end{tabular}}
\newcommand{\bdoc}{\begin{document}}
\newcommand{\edoc}{\end{document}}
\newcommand{\beqy}{\begin{eqnarray}}
\newcommand{\eeqy}{\end{eqnarray}}

\newcommand{\beqi}{\begin{eqnarray*}}
\newcommand{\eeqi}{\end{eqnarray*}}
\newcommand{\bitem}{\begin{itemize}}

\newcommand{\eitem}{\end{itemize}}
\newcommand{\nln}{\newline}
\newcommand{\newt}{\newtheorem}

\newcommand{\pa}{\partial}
\newcommand{\re}{{I\!\!R}}
\newcommand{\ren}{\re^N}
\newcommand{\xr}{x\in\re }
\newcommand{\x}{\times}
\newcommand{\dyle}{\displaystyle}
\newcommand{\ene}{{I\!\!N}}
\newcommand{\irn}{\int\limits_{\re^N}}
\newcommand{\io}{\int\limits_{\O}}
\newcommand{\meas}{{\rm meas\,}}

\newcommand{\tx}[1]{\mbox{\;{#1}\;}}
\newcommand{\erre}{\mathbb{R}}
\newcommand{\Hess}{\operatorname{Hess}}
\newcommand{\metric}{g}
\newcommand{\sign}{{\rm sign}}
\newcommand{\map}{\longrightarrow }
\newcommand{\imp}{\Longrightarrow }
\renewcommand{\div}{\nabla\cdot }
\newcommand{\sen}{{\rm sen\,}}
\newcommand{\tg}{{\rm tg\,}}
\newcommand{\arcsen}{{\rm arcsen\,}}
\newcommand{\arctg}{{\rm arctg\,}}
\newcommand{\supp}{{\textsl supp\ }}
\newcommand{\ity}{\int_{-\iy}^{+\iy}}
\newcommand{\limit}{\lim\limits}
\newcommand{\limi}{\limit_{n\to\infty}}
\newcommand{\sumi}{\sum\limits_{n=1}^{\infty}}
\newcommand{\ulu}{\underline u}
\newcommand{\ulw}{\underline w}
\newcommand{\ulz}{\underline z}
\newcommand{\ulv}{\underline v}
\newcommand{\uls}{\underline s}
\newcommand{\olu}{\overline u}
\newcommand{\olv}{\overline v}
\newcommand{\ols}{\overline s}
\newcommand{\ob}{\overline\b}
\newcommand{\ovar}{\overline\var}
\newcommand{\wv}{\widetilde v}
\newcommand{\wu}{\widetilde u}
\newcommand{\ws}{\widetilde s}
\renewcommand{\a }{\alpha }
\renewcommand{\b }{\beta }
\newcommand{\g }{\gamma}
\newcommand{\G }{\Gamma }
\renewcommand{\d }{\delta }

\newcommand{\loc}{{\rm loc}}
\newcommand{\D }{\Delta }
\newcommand{\e }{\varepsilon }
\newcommand{\z }{\zeta }
\renewcommand{\l }{\lambda }
\renewcommand{\L }{\Lambda }
\newcommand{\m }{\mu }
\newcommand{\n }{\nabla }
\newcommand{\s }{\sigma }
\newcommand{\Sig }{\Sigma }
\renewcommand{\t }{\tau }
\newcommand{\var }{\varphi }
\renewcommand{\o }{\omega }
\renewcommand{\O }{\Omega }
\newcommand{\bR}{{\bf R}}
\newcommand{\bC}{{\bf C}}
\newcommand{\bZ}{{\bf Z}}
\newcommand{\bN}{{\bf N}}
\newcommand{\bQ}{{\bf Q}}
\newcommand{\bK}{{\bf K}}
\newcommand{\bI}{{\bf I}}
\newcommand{\bv}{{\bf v}}
\newcommand{\bV}{{\bf V}}

\newcommand{\pair}[1]{g\left(#1\right)}

\def\qed{\unskip\kern 6pt \penalty 500
\raise -2pt\hbox{\vrule \vbox to10pt{\hrule width 4pt
\vfill\hrule}\vrule}\par}

\newenvironment{Proof}{\removelastskip\vskip12pt
plus 1pt \noindent\em\rm}{\hfill {\qed \hskip .2cm}}

\title[Hyperbolic problems on riemannian manifolds]{Nonexistence for hyperbolic problems \\ on riemannian manifolds}
\author{Dario D. Monticelli}
\author{Fabio Punzo}
\author{Marco Squassina}

\address[D.\ Monticelli]{Dipartimento di Matematica \newline\indent
	Politecnico di Milano, Milano, Italy}
\email{dario.monticelli@polimi.it}

\address[D.\ Monticelli]{Dipartimento di Matematica \newline\indent
	Politecnico di Milano, Milano, Italy}
\email{fabio.punzo@polimi.it}

\address[M.\ Squassina]{Dipartimento di Matematica e Fisica \newline\indent
	Universit\`a Cattolica del Sacro Cuore \newline\indent
	Via dei Musei 41, I-25121 Brescia, Italy}
\email{marco.squassina@unicatt.it}


\thanks{The authors are members of the Gruppo Nazionale per l'Analisi Matematica, la Probabilit\`a e le loro Applicazioni.}

\begin{document}

\begin{abstract}
We establish necessary conditions for the existence of solutions to a class of semilinear hyperbolic
problems defined on complete noncompact Riemannian manifolds, extending some nonexistence results
for the wave operator with power nonlinearity on the whole euclidean space.  A general weight
function depending on spacetime is allowed in front of the power nonlinearity.
\end{abstract}

\maketitle

\begin{center}
	\begin{minipage}{8.5cm}
		\small
		\tableofcontents
	\end{minipage}
\end{center}

\medskip

\section{Introduction} \setcounter{equation}{0}
In the framework of parabolic equations, a classical result obtained by H.\ Fujita
\cite{fuji} in 1966 (see also \cite{haya}) states that the semilinear problem on $\mathbb R^n$
\begin{equation}\label{eq1}
\left\{
\begin{array}{ll}
\,  u_{t} - \Delta u\, =u^q \, &\textrm{in}\,\, \mathbb R^n\times (0,\infty),
\\ & \\
\textrm{ }u \, = u_0& \textrm{in}\,\,  \mathbb R^n\times \{0\},
\end{array}
\right.
\end{equation}
does not admit nontrivial global solutions provided that
$$
1<q\leq 1+\frac{2}{n},\qquad
u_0\in L^1(\mathbb R^n),\quad u_0\geq 0,\quad  u_0\not \equiv 0.
$$
This range of values of the power $q$ is often referred to as the {\it blow-up range} and
$1+{2/n}$ is called the {\it Fujita exponent}. Many generalization of this result were
derived after Fujita's paper, see for instance the 2000 survey paper by K.\ Deng and H.A.\ Levine \cite{survey}
and the references therein. More recently, an extension of this kind of achievements when the Euclidean
space is replaced by noncompact Riemannian manifolds was obtained in \cite{BPT, MMP2, Pu1, zhang}, under suitable geometric hypotheses.
Similar results for elliptic equations on noncompact Riemannian manifolds have been also investigated (see, e.g., \cite{GrigKond}, \cite{GrigS}, \cite{MMP1}).

As far as {\em hyperbolic} equations are concerned, the situation is quite
different compared to the parabolic setting
and, to our knowledge, the first contribution on the subject of nonexistence
for the wave equation on $\mathbb R^n$ is \cite{kato, shae} and more generally \cite{MP}, where it is was proved
that the problem
\begin{equation}\label{eq-mod}
\left\{
\begin{array}{ll}
\,  u_{tt} - \Delta u\, \geq \,  |u|^q \, &\textrm{in}\,\, \mathbb R^n\times (0,\infty),
\\&\\
\textrm{ }u \, = u_0& \textrm{in}\,\,  \mathbb R^n\times \{0\},\\&\\
\textrm{ }u_t\, = u_1&\textrm{in}\,\,  \mathbb R^n\times \{0\},
\end{array}
\right.
\end{equation}
admits no nontrivial solution, provided that
$$
1<q\leq \frac{n+1}{n-1},\qquad \liminf_{R\to\infty}\int_{B_R(0)} u_1 dx \geq 0\,.
$$
Therefore, in some sense, the exponent $(n+1)/(n-1)$ plays in the hyperbolic case the same
r\v ole played by the Fujita exponent $1+{2/n}$ in the parabolic case.

\smallskip

The main goal of this paper is to obtain the counterpart of the result in
\cite{MP} on complete noncompact Riemannian manifolds, as it was done in \cite{MMP2} for the parabolic case from the original Fujita breakthrough
achievement. More precisely, we consider the hyperbolic problem
\begin{equation}\label{eq1}
\left\{
\begin{array}{ll}
 \,  u_{tt} - \Delta u\, \geq \, V |u|^q \, &\textrm{in}\,\,M\times (0,\infty),
\\&\\
\textrm{ }u \, = u_0& \textrm{in}\,\,  M\times \{0\},\\&\\
\textrm{ }u_t\, = u_1&\textrm{in}\,\,  M\times \{0\}.
\end{array}
\right.
\end{equation}
where $M$ is a complete noncompact
Riemannian manifold of dimension $n$, endowed with a metric tensor $g$, $\Delta$ is the Laplace-Beltrami operator on $M$,
$q>1$, $V\in L^1_\textrm{loc}(M\times [0, \infty))$ and $V>0$ almost everywhere.

For a fixed reference point $o\in M$, we set $r(x) :=\operatorname{dist}(x, o)$, and
for any $x\in M$ and for any $R>0$ we set $B_R(o):=
\big\{x\in M\,:\, r(x)<R\,\big\}$.

In \cite{Ru} some nonexistence results for problem \eqref{eq1} have been stated; however, in the proofs in \cite{Ru} it is implicitly assumed that $u\geq 0$ in $M\times (0, \infty)$ and that $x\mapsto r(x)$ is of class $C^2$ in $M\setminus\{o\}$. Observe that the hypothesis that $u\geq 0$ in $M\times (0, \infty)$ is not natural for solutions of hyperbolic equations; in fact, also in \cite{MP}, where $M=\mathbb R^n$, sign-changing solutions were considered.
Moreover, the required regularity for the function $x\mapsto r(x)$ on general Riemannian manifolds is not guaranteed; in particular, it holds whenever the cut locus of the point $o$ is empty. In this paper, on the one hand, we do not require that the function $x\mapsto r(x)$ is of class $C^2$ in $M\setminus\{o\}$ (see Theorem \ref{teo1} below). On the other hand, we remove the unnatural condition that $u\geq 0$. In addition, we assume weaker hypotheses both on the volume growth of geodesic balls and on the potential.

We establish nonexistence of very weak solutions of problem \eqref{eq1}, under a suitable bound from below on the Ricci curvature and an appropriate weighted volume growth condition on geodesic balls, with weight depending on the potential. The methods of proofs that we use present some important differences with respect to the elliptic and the parabolic cases. Indeed,
in \cite{GrigKond}, \cite{GrigS}, \cite{MMP1}, \cite{MMP2} local weak {\it nonnegative} solutions have been considered, and the argument that yields nonexistence of solutions was based on a careful choice of radial test functions, i.e. depending on $r(x)$. Observe that only the gradient of such test functions was involved, but not their second derivatives; indeed, on general manifolds such a gradient is well-defined, since $|\nabla r(x)|=1$ a.e. in $M$. In the present situation those techniques do not work, loosely speaking, due to the presence of the term $u_{tt}$ and to the fact that $u$ can change sign. Hence, one needs to consider {\em very weak} solutions of problem \eqref{eq1} (see Definition \ref{defsol} below), and consequently one must estimate also second derivatives of the test functions. This prevents us to use test functions depending on $r(x)$, since in general $x\mapsto r(x)$ is not $C^2$ in $M\setminus\{o\}$. In order to overcome this obstacle, we shall use appropriate regular test functions introduced very recently in the nice paper \cite{BS}. We also address separately Cartan-Hamard maninfolds. In this special case, since $\operatorname{Cut}(o)=\emptyset$, we directly use test functions depending on $r(x)$; furthermore, we can allow the Ricci curvature to diverge negatively faster than in the general case.

\smallskip
The paper is organized as follows. In Section \ref{mr} we state our main results, we compare them with those known in $\mathbb R^n$, and we discuss some examples; Section \ref{RG} is devoted to some notions of Riemannian Geometry that will be used in the sequel. In Section \ref{prel} we give the precise definition of solutions to problem \eqref{eq1}, we recall an important useful result from \cite{BS}, then we obtain key a priori estimates for solutions of problem \eqref{eq1}. Finally, we prove the main results in Section \ref{proof}.

\section{Main results}\label{mr}

\subsection{Case $\boldsymbol{\sigma\geq -2}$}
The following is one of the main result of the paper.

\begin{theorem}[Case $\boldsymbol{\sigma\geq -2}$]
	\label{teo1}
Let $q>1$, $V\in L^1_{\loc}(M\times [0, \infty))$, $V>0$ a.e. in $M\times (0, \infty)$, $u_0\in L^1_{\loc}(M)$, $u_1\in L^1_{\loc}(M)$ and
\begin{equation}\label{eq5}
\int_{M} u_1 d\mu \geq 0\,.
\end{equation}
Suppose that, for some $C_0\geq 0$ and $\sigma\in [-2, \infty)$,
\begin{equation}\label{eq4}
\operatorname{Ric}\geq - (n-1)\frac{C_0}{(1+r(x)^2)^{\frac{\sigma}2}}\,.
\end{equation}
Assume that there exist $R_0>0$ and $\bar C>0$ such that, for every $R>R_0$,
\begin{equation}\label{eq11}
\int_0^{R}\int_{B_R(o)} V^{-\frac{1}{q-1}}\, d\mu dt \leq \bar C R^{\min\left\{1+\frac{\sigma}2,2\right\}\frac{q}{q-1}}\,.
\end{equation}
Then the unique solution of problem \eqref{eq1} is $u\equiv 0$.
\end{theorem}

Compared to the Euclidean case $M={\mathbb R}^n$, in the general case the main problem is that of detecting suitable
well-defined test functions. To this aim, we exploit a recent result from \cite{BS}, stating that
under condition \eqref{eq5} there exists a family  $\{\phi_R\}_{R\geq 1}\subset C^{\infty}_c(M)$
with $0\leq \phi_R\leq 1$ in $B_R(o)$, $\phi_R\equiv 1$ in $B_R(o)$, $\operatorname{supp} \phi_R\subset B_{\gamma R}(o)$ for all $\gamma>\gamma_0$,
with $\gamma_0>1$ independent of $R$ and
$$
|\nabla\phi_R|\leq \frac{C}R,\qquad
|\Delta \phi_R|\leq \frac{C}{R^{1+\frac{\sigma}{2}}},
$$
for some constant $C$ independent of $R$. In fact, basically, in the proof
we test the equation by $\phi_R(x)\eta_R(t),$ where $\eta_R\in C^{2}([0, \infty))$
is a suitable cut-off time function and apply a technical result (see Lemma \ref{lemma2})
to the compact set $K_R=[0, R]\times\overline{B_R(o)}$, which finally leads to
\begin{equation*}
\int_{K_R} |u|^q V\,d\mu dt\leq - \int_Mu_1(x,0)\psi^\alpha(x,0)\,d\mu + C \left(\int_{M\times (0, \infty)\setminus K_R} |u|^qV\,d\mu dt \right)^{\frac 1 q}\,,
\end{equation*}
for some $C>0$ independent of $R$. This finally yields the assertion $u=0$ by letting $R\to\infty$.

\bigskip
\noindent
We notice that condition \eqref{eq5} can be weakened to
\[
\limsup_{R\to \infty} \left\{\int_{B_R(o)} u_1^+ d\mu - \int_{B_{\gamma R}} u_1^- d\mu \right\} \geq 0\,,
\]
as it can be seen from the proof. Here $\gamma>0$ is sufficiently large (cf.\ Proposition \ref{propSB} below).

%

If the potential $V$ is written as a product of two functions, one only depending on the space variables and the other only on the time variable, then by our result we can deduce the following consequences.

\begin{corollary}[Splitting case]
	\label{cor2}
Let $q>1$, $\sigma\in [-2, 2]$, $V\in L^1_{\loc}(M\times [0, \infty))$,
$V>0$ a.e. in $M\times (0, \infty)$, $u_0\in L^1_{\loc}(M)$, $u_1\in L^1_{\loc}(M)$. Assume that \eqref{eq4} holds, and that
\[V(x,t)=f(x)g(t),\qquad x\in M,\, t>0\,. \]
Moreover, assume that there exist $R_0>0$ and $\bar C>0$ such that, for every $R>R_0$,
\[ \int_0^R g(t)^{-\frac{1}{q-1}}dt \leq \bar C R^{\alpha}\,,\quad  \int_{B_R}f(x)^{-\frac{1}{q-1}}\, d\mu\leq \bar C R^{\beta}\,, \]
with $\alpha, \beta\in \mathbb R\,, \alpha+\beta\leq \left(1+\frac{\sigma}{2}\right){\frac{q}{q-1}}$\,.
Then the unique solution of problem \eqref{eq1} is $u\equiv 0$.
\end{corollary}

\begin{remark}[Model case]
	\rm
	\label{oss1}
	By Corollary \ref{cor2} we can obtain a nonexistence result for problem \eqref{eq1} with $V\equiv 1$.
Assume that \eqref{eq4} holds, and that there exist $R_0>0$ and $\bar C>0$ such that, for every $R>R_0$,
\begin{equation}\label{eq30}
\mu(B_R(o))\leq \bar C R^{\frac 1{q-1}\min\left\{1+\frac{\sigma q}{2},\, q+1\right\}}.
\end{equation}
Then, by Corollary \ref{cor2}, problem \eqref{eq1} with $V\equiv 1$ has only the trivial solution $u\equiv 0$.
Furthermore, when $M=\mathbb R^n$, \eqref{eq4} is satisfied with $\sigma=2$. Since in this case
$\mu(B_R(o))=\bar C R^n$, condition \eqref{eq30} holds, whenever
\begin{equation}\label{eq31}
1<q\leq \frac{n+1}{n-1}\,.
\end{equation}
Thus we recover the result established in \cite{MP}.
\end{remark}

\begin{remark}[Examples of $V$]
	\rm
	\label{oss2}
	By using the co-area formula and Bishop-Gromov volume comparison theorems, one can prove the following sufficient conditions for hypothesis \eqref{eq11}.
\begin{itemize}
\item[$(i)$] Assume that \eqref{eq4} holds for some $\sigma\in [-2, 2)$, and that, for some $C>0$,
\[ V(x, t)\geq C r(x)^{\beta}\exp\{\alpha r(x)^{1-\frac{\sigma}2}\}\,, \quad x\in M, t>0\,.\]
Then, by volume comparison theorems (see e.g. \cite{Grig}),
\[\mu(B_R(o))\leq \bar C \exp\{B(n-1) R^{1-\frac{\sigma}2}\}\,,\]
for some $B>0$ and $\bar C>0$. Hence, if
\[\alpha\geq B(n-1)(q-1),\quad \beta\geq q\left(1-\frac{\sigma}2\right)-2\,,\]
then condition \eqref{eq11} is satisfied, and thus Theorem \ref{teo1} applies.

\item[$(ii)$] Now, assume that \eqref{eq4} holds for $\sigma=2$, and that, for some $C>0$,
\[ V(x, t)\geq C r(x)^{\beta}\,, \quad x\in M, t>0\,.\]
Then, by volume comparison theorems (see e.g. \cite{Grig}),
\[\mu(B_R(o))\leq \bar C R^{(n-1)\delta+1} \,,\]
with $\delta=\frac{1+\sqrt{1+\frac{4C_0}{n-1}}}{2}$, for some $\bar C>0$. Hence, if
\[\beta\geq (n-1)\delta(q-1)-2\,,\]
then condition \eqref{eq11} is satisfied, and thus Theorem \ref{teo1} applies.

\item[$(iii)$] Suppose that \eqref{eq4} holds for some $\sigma>2$, and that, for some $C>0$,
\[ V(x, t)\geq C r(x)^{\beta}\,, \quad x\in M, t>0\,.\]
Then, by volume comparison theorems (see e.g. \cite{Grig}),
\[\mu(B_R(o))\leq \bar C R^{n} \,,\]
for some $\bar C>0.$ Hence, if
\[\beta\geq (n-1)(q-1)-2\,,\]
then condition \eqref{eq11} is satisfied, and thus Theorem \ref{teo1} applies.

\item[$(iv)$] If
\[\operatorname{Ric}\geq 0\,,\]
then the same result as in $(iii)$ holds. Then, this case is equivalent to the choice $M=\mathbb{R}^n$.

\item[$(v)$] If in $(iii)$ or in $(iv)$ we further assume that $\beta=0$, then \eqref{eq11} is satisfied, provided that \eqref{eq31} holds.
\end{itemize}
\end{remark}

\subsection{Case $\boldsymbol{\sigma<-2}$}
Now we consider Cartan-Hadamard manifolds (see subsection \ref{CH} below). In this special framework we can also consider the case that
\eqref{eq4} is satisfied with $\sigma<-2$.
 \begin{theorem}[Case $\boldsymbol{\sigma<-2}$]
 	\label{teo2} Let $M$ be a Cartan-Hadamard manifold.
Let $q>1$, $V\in L^1_{\loc}(M\times [0, \infty))$, $V>0$ a.e. in $M\times (0, \infty)$, $u_0\in L^1_{\loc}(M)$, $u_1\in L^1_{\loc}(M)$, in addition assume that \eqref{eq5} holds. Suppose that there exist $C_0>0, \sigma<-2$ such that
\begin{equation}\label{eq34}
\operatorname{Ric}\left(\frac{\partial}{\partial r}, \frac{\partial}{\partial r}\right) \geq - (n-1)\frac{C_0}{(1+r(x))^{\sigma}}\,.
\end{equation}
Assume that there exist $R_0>0$ and $\bar C>0$ such that, for every $R>R_0$,
\begin{equation}\label{eq35}
\int_0^{R}\int_{B_R(o)} V^{-\frac{1}{q-1}}\, d\mu dt \leq \bar C R^{\left(1+\frac{\sigma}2\right)\frac{q}{q-1}}\,.
\end{equation}
Then the unique solution of problem \eqref{eq1} is $u\equiv 0$.
\end{theorem}

\begin{remark}[Examples of $V$]
	\rm
Let $M$ be a Cartan-Hadamard manifold, suppose that \eqref{eq34} holds for some $\sigma<-2$. Assume that for some $C>0$,
\[ V(x, t)\geq C r(x)^{\beta}\exp\{\alpha r(x)^{1-\frac{\sigma}2}\}\,, \quad x\in M, t>0\,.\]
Then, arguing as in Remark \ref{oss2}-(i), one has
\[\mu(B_R(o))\leq \bar C \exp\{B(n-1) R^{1-\frac{\sigma}2}\}\,,\]
for some $B>0$ and $\bar C>0$. So, if
\[\alpha\geq B(n-1)(q-1),\quad \beta\geq q\left(1-\frac{\sigma}2\right)-2\,,\]
then condition \eqref{eq35} is satisfied, and thus Theorem \ref{teo2} applies.
\end{remark}

\section{Basics from Riemannian geometry}\label{RG}
For the reader's convenience we first recall  some notions and results from Riemannian geometry, see e.g. \cite{AMR}. Let $M$ be a Riemannian manifold of dimension $m$ endowed with a metric $\metric=\langle\cdot,\cdot\rangle$. Let $d\mu$ be the canonical Riemannian measure on $M$.  We denote by $x$ an arbitrary point of $M$ and
let $x^1,\ldots, x^m$ be the coordinate functions in  the local chart $U$. Then we have
\begin{equation}
\label{GP1.1}
\metric=g_{ij}\,dx^{i}\otimes dx^{j}
\end{equation}
where $dx^{i}$ denotes the differential of the function $x^{i}$ and
$g_{ij}$ are the (local) components of the metric, defined by
$g_{ij}=\langle\frac{\partial}{\partial x^{i}},
\frac{\partial}{\partial x^{j}}\rangle$.
We will denote by $[g^{ij}]$ the inverse of the matrix $[g_{ij}]$. In the sequel we
shall use the Einstein summation convention over repeated indices.
For a fixed point $o\in M$, we denote $r(x) :=\operatorname{dist}(x, o)$
for any $x\in M$ and for any $R>0$ we set
$$
B_R(o):=
\big\{x\in M\,:\, r(x)<R\,\big\}.
$$
For any smooth function $u:M\to \erre$, the \emph{gradient} of $u$ relative to the metric $g$ of $M$, $\nabla u$, is the vector field dual to the $1$-form $du$, that is
\[
\langle\nabla u, X\rangle = du(X)=X(u)
\]
for all smooth vector fields $X$ on $M$. In local coordinates
we have
$$
\nabla u=u^i\frac{\partial}{\partial x^i} \tx{with}
u^j=g^{ij}\frac{\partial u}{\partial x^i},
$$
and
$$|\nabla u|^2=\langle\nabla u, \nabla u\rangle=g^{ij}\frac{\partial u}{\partial x_i}\frac{\partial u}{\partial x_j}.$$
The \emph{divergence} of a vector field $X$ on $M$ is given by the trace of $\nabla X$, the covariant derivative of $X$, where $\nabla$ is the (unique) Levi-Civita connection associated to the metric $\metric$. If $X=X^i\frac{\partial}{\partial x^i}$,
the divergence of $X$ can be expressed in local coordinates as
$$
\operatorname{div}X= \frac{\partial X^i}{\partial
	x^i}+X^k\Gamma^i_{ki},
$$
where $\Gamma^k_{ij}$ are the Christoffel symbols
\[
\Gamma^k_{ij} =
\frac{1}{2}g^{kl}\left(\frac{\partial g_{il}}{\partial x^j}+\frac{\partial g_{jl}}{\partial x^i}-\frac{\partial g_{ij}}{\partial x^l}\right).
\]

The Hessian of $u$  is defined as the $2$--tensor $\Hess(u)=\nabla du$, the covariant derivative of $du$, and its components $u_{ij}$ are in local coordinates

$$
u_{ij}=\frac{\partial ^2u}{\partial x^i\partial x^j}
-\Gamma^k_{ij}\frac{\partial u}{\partial x^k}.
$$

The
\emph{Laplace--Beltrami operator} of $u$ is the trace of the Hessian, or equivalently the divergence of the gradient, i.e.
\[
\label{GP3.4}
\Delta u =\operatorname{Tr}(\Hess(u))=\operatorname{div}(\nabla u).
\]
In local coordinates it has the form
$$
\Delta u= \mathfrak g^{-1}\frac{\partial}{\partial
	x^i}\left(\mathfrak g \, g^{ij} \frac{\partial u}{\partial
	x^j}\right),\quad \tx{where} \mathfrak{g}=\sqrt{{\rm det}(g_{ij})}.
$$

We denote by  $\operatorname{Ric}$ the \emph{Ricci tensor} which  is expressed in local coordinates as
\[
R_{ij}=R_{ji}=\frac{\partial\Gamma^l_{ij}}{\partial x^l}-\frac{\partial\Gamma^{l}_{il}}{\partial x^j}+\Gamma^k_{ij}\Gamma^l_{kl}-\Gamma^k_{il}\Gamma^l_{kj},
\]
and we write $\operatorname{Ric}\geq h(x)$ for a given function $h:M\rightarrow\mathbb{R}$ to intend
$$
\operatorname{Ric}(X,X)\geq h(x)|X|^2,\qquad\textrm{for every vector field }X.
$$
Let
$$\mathcal A:=\left\{f\in C^\infty((0,\infty))\cap C^1([0,\infty)): \, f'(0)=1, \, f(0)=0, \, f>0 \ \textrm{in}\;\, (0,\infty)\right\} .$$
If
\[\operatorname{Ric}\left(\frac{\partial}{\partial r}, \frac{\partial}{\partial r}\right)\geq - (n-1)\frac{\varphi''(r(x))}{\varphi(r(x))}\quad \textrm{for any}\;\; x\in M\setminus\{o\}\,,\]
for some $\varphi\in \mathcal A$, then, by volume comparison theorem,
\begin{equation}\label{eq32}
\operatorname{Vol}(B_R(o)) \leq  C\int_0^R \varphi^{n-1}(\xi)\, d\xi\,,
\end{equation}
for some $C>0$ independent of $R$.

\subsection{Cartan-Hadamard manifolds}\label{CH} In the sequel we consider {\it Cartan-Hadamard} manifolds, i.e. simply connected complete noncompact Riemannian manifolds with
nonpositive sectional curvatures. Observe that (see, e.g. \cite{Grig}) for a Cartan-Hadamard manifold $M$
the {\it cut locus } of $o$, $\operatorname{Cut}(o)$, is empty for any point $o \in M,$ thus $M$ is a manifold with a pole.
For any $x\in M\setminus \{o\}$, one can define the {\it polar
coordinates} with respect to $o$. Namely, for
any point $x\in M\setminus\{o\}$
there correspond a polar radius $r(x) :=\operatorname{dist}(x, o)$
and a polar angle $\theta\in \mathbb S^{n-1}$ such that the shortest
geodesics from $o$ to $x$ starts at $o$ with direction $\theta$ in
the tangent space $T_o M$. Since we can identify $T_o M$ with
$\mathbb R^n$, $\theta$ can be regarded as a point of $\mathbb
S^{n-1}.$

The Riemannian metric in $M\setminus\{o\}$ in polar coordinates reads
\[g= dr^2+A_{ij}(r, \theta)d\theta^i d\theta^j, \]
where $(\theta^1, \ldots, \theta^{n-1})$ are coordinates in $\mathbb S^{n-1}$ and $(A_{ij})$ is a positive definite matrix.
It is not difficult to see that the Laplace-Beltrami operator in polar coordinates has the form
\begin{equation}\label{eq36} \Delta = \frac{\partial^2}{\partial r^2} +
\mathcal F(r, \theta)\frac{\partial}{\partial r}+\Delta_{S_{r}},
\end{equation}
where $\mathcal F(r, \theta):=\frac{\partial}{\partial
r}\big(\log\sqrt{A(r,\theta)}\big)$, $A(r,\theta):=\det
(A_{ij}(r,\theta))$, $\Delta_{S_r}$ is the Laplace-Beltrami operator
on the submanifold $S_{r}:=\partial B_r(o)$\,.

A manifold with a pole is a {\it spherically symmetric manifold} or a {\it model}, if the Riemannian metric is given by
\begin{equation}\label{eq37}
g= dr^2+ f^2(r)d\theta^2,
\end{equation}
where $d\theta^2=\beta_{ij}d\theta^i d \theta^j$ is the standard
metric in $\mathbb S^{n-1}$, $\beta_{ij}$ being smooth functions of
$\theta^1, \ldots, \theta^{n-1},$ and $f\in \mathcal A$.
In this case, we write $M\equiv
M_f$; furthermore, we have $\sqrt{A(r,\theta)}=f^{n-1}(r)$, so
that
\begin{equation}\label{eq38}
\Delta = \frac{\partial^2}{\partial r^2}+ (n-1)\frac{f'}{f}\frac{\partial}{\partial r}+ \frac1{f^2}\Delta_{\mathbb S^{n-1}}\,,
\end{equation}
where $\Delta_{\mathbb S^{n-1}}$ is the Laplace-Beltrami operator in
$\mathbb S^{n-1}\,.$ Observe that for $f(r)=r$, $M=\mathbb R^n$, while for $f(r)=\sinh r$, $M$ is the $n-$dimensional hyperbolic space $\mathbb H^n$.

\smallskip

Let us recall comparison results for Ricci curvature
that will be used in the sequel. Observe that (see
\cite[Section 15]{Grig}), if
\begin{equation}\label{eq41}
\textrm{Ric}\left(\frac{\partial}{\partial r}, \frac{\partial}{\partial r}\right)\geq -(n-1)\frac{\phi''(r)}{\phi(r)}\quad \textrm{for all}\;\; x=(r,\theta)\in M\setminus\{o\},
\end{equation} for some function $\phi\in \mathcal A$, then
\begin{equation}\label{eq42}
\mathcal F(r, \theta)\leq (n-1)\frac{\phi'(r)}{ \phi(r)}\quad \textrm{for all}\;\; r>0, \theta \in \mathbb S^{n-1}\,\,\textrm{with}\,\, x=(r,
\theta)\in M\setminus \{o\}\,.
\end{equation}
On the other hand, since $M$ is a Cartan-Hadamard manifold, $\mathcal F(r, \theta)\geq 0$.
If $M_f$ is a model manifold, then for any $x=(r,
\theta)\in M_f\setminus\{o\}$
\[\textrm{Ric}\left(\frac{\partial}{\partial r}, \frac{\partial}{\partial r}\right)=-(n-1)\frac{f''(r(x))}{f(r(x))}\,.\]

\section{Preliminaries and a priori estimates}\label{prel}

\begin{definition}\label{defsol}
We say that $u\in L^1_\textrm{loc}(M\times [0, +\infty))\cap L^q_\textrm{loc}(M\times [0, +\infty),Vd\mu)$ is a (very weak) solution of problem \eqref{eq1} if
\begin{eqnarray}
&&\label{eq2} \int_0^{+\infty}\!\int_M\varphi|u|^qV\,d\mu dt\leq \\
&&\nonumber\qquad \int_0^{+\infty}\!\int_Mu(\varphi_{tt}-\Delta\varphi)\,d\mu dt-\int_Mu_1(x)\varphi(x,0)\,d\mu+\int_Mu_0\varphi_t(x,0)\,d\mu
\end{eqnarray}
for every $C^2$ nonnegative function $\varphi$ with compact support in $M\times [0, +\infty)$.
\end{definition}

We use the following result established in \cite{BS}.
\begin{proposition}\label{propSB}
Suppose that condition \eqref{eq4} is satisfied with $\sigma \in [-2, 2]$. Then there exists a family of functions $\{\phi_R\}_{R\geq 1}\subset C^{\infty}_c(M)$ with the following properties:
\begin{itemize}
\item[$(i)$] $0\leq \phi_R\leq 1$ in $M$ and $\phi\equiv 1$ in $B_R(o)$;
\item[$(ii)$] $\operatorname{supp} \phi_R\subset B_{\gamma R}(o)$ for all $\gamma>\gamma_0$, for some $\gamma_0>1$ independent of $R$;
\item[$(iii)$] $|\nabla\phi_R|\leq \frac{C}R$ for some constant $C$ independent of $R$;
\item[$(iv)$] $|\Delta \phi_R|\leq \frac{C}{R^{1+\frac{\sigma}{2}}}$ for some constant $C$ independent of $R$.
\end{itemize}
\end{proposition}

\medskip

In the following two lemmas we obtain two a priori estimates for solutions of problem \eqref{eq1} that will play a crucial role in the proofs of Theorems \ref{teo1}, \ref{teo2}.

\begin{lemma}\label{lemma1}
  Let  $u\in L^1_{loc}(M\times [0, +\infty))\cap L^q_\textrm{loc}(M\times [0, +\infty),Vd\mu)$ be a (very weak) solution of problem \eqref{eq1}, $\alpha>\frac{2q}{q-1}$, $\psi\in C^2(M\times[0,\infty))$ nonnegative and compactly supported. Then
  \begin{equation}\label{eq3}
  \begin{aligned}
&\int_0^\infty\!\int_M |u|^q\psi^\alpha V\,d\mu dt\leq -\frac{q}{q-1}\int_Mu_1(x,0)\psi^\alpha(x,0)-u_0(x,0)(\psi^\alpha(x,0))_t\,d\mu \\ &+\int_0^\infty\!\int_MV^{-\frac{1}{q-1}}\psi^{\alpha-\frac{2q}{q-1}}\big|\alpha(\alpha-1)\psi_t^2+\alpha\psi\psi_{tt}-\alpha(\alpha-1)|\nabla\psi|^2-\alpha\psi\Delta\psi\big|^\frac{q}{q-1}\,d\mu dt
  \end{aligned}
  \end{equation}
\end{lemma}

\begin{proof}
  We use the definition of solution of problem \eqref{eq1}, using the admissible test function $\varphi=\psi^{\alpha}$. We get
\[\int_0^{\infty}\int_M |u|^q\psi^{\alpha} V\, d\mu dt \leq  \int_0^{\infty}\int_M |u||(\psi^{\alpha})_{tt}-\Delta(\psi^{\alpha})| d\mu dt
\]
\[+\int_M [u_0(x) (\psi^{\alpha})_t(x,0)-u_1(x)\psi^{\alpha}(x,0)]\,d\mu
\]
\[\leq \frac{1}{q} \int_0^{\infty}\int_M |u|^q \psi^{\alpha} V\, d\mu dt + \frac{q}{q-1} \int_0^{\infty}\int_M|(\psi^{\alpha})_{tt} - \Delta(\psi^{\alpha})|^{\frac{q}{q-1}}(V \psi^{\alpha})^{-\frac{1}{q-1}}\, d\mu dt\]
\[ + \int_M [u_0(x) (\psi^{\alpha})_t(x,0)- u_1(x)\psi^{\alpha}(x,0)]\,d\mu \,.
\]
Then
\[\int_0^{\infty}\int_M |u|^q \psi^{\alpha} V\, d\mu dt \leq \int_0^{\infty}\int_M|(\psi^{\alpha})_{tt} - \Delta(\psi^{\alpha})|^{\frac{q}{q-1}}(V \psi^{\alpha})^{-\frac{1}{q-1}}\, d\mu dt\]\[+
\frac{q}{q-1} \int_M [u_0(x) (\psi^{\alpha})_t(x,0)- u_1(x)\psi^{\alpha}(x,0)]\,d\mu\,.
\]
Observe that
\[(\psi^{\alpha})_{tt}=\alpha(\alpha-1)\psi^{\alpha-2}(\psi_t)^2 + \alpha \psi^{\alpha-1}\psi_{tt}\,;\]
\[\Delta(\psi^{\alpha})=\alpha(\alpha-1)\psi^{\alpha-2}|\nabla\psi|^2 + \alpha \psi^{\alpha-1}\Delta\psi\,.\]
Thus we have
\[  \int_0^{\infty}\int_M |u|^q \psi^{\alpha} V\, d\mu dt \]\[\leq \int_0^{\infty}\int_M V^{-\frac{1}{q-1}}\psi^{\alpha-\frac{2q}{q-1}}\left|\alpha(\alpha-1)(\psi_t)^2+\alpha\psi\psi_{tt}-\alpha(\alpha-1)|\nabla\psi|^2-\alpha\psi\Delta\psi \right|^{\frac{q}{q-1}}\, d\mu dt\]
\[+\frac{q}{q-1}\int_M [u_0(x) (\psi^{\alpha})_t(x,0)- u_1(x)\psi^{\alpha}(x,0)]\,d\mu\,.
\]
This concludes the proof.
\end{proof}

\begin{lemma}\label{lemma2}
Let  $u\in L^1_{loc}(M\times [0, +\infty))\cap L^q_\textrm{loc}(M\times [0, +\infty),Vd\mu)$ be a (very weak) solution of problem \eqref{eq1}, $\alpha>\frac{2q}{q-1}$. Let $K\subset M\times [0, \infty)$ be a compact subset, $\psi\in C^2(M\times[0,\infty))$ nonnegative and compactly supported, $\psi\equiv 1$ in $K$. Set $S_K:=[M\times (0, \infty)]\setminus K$. Then
  \begin{equation}\label{eq10}
\begin{aligned}
&\int_0^\infty\!\int_M |u|^q\psi^\alpha V\,d\mu dt\leq - \int_Mu_1(x,0)\psi^\alpha(x,0)-u_0(x,0)(\psi^\alpha(x,0))_t\,d\mu \\ &+\left(\int_{S_K} V^{-\frac{1}{q-1}}\psi^{\alpha-\frac{2q}{q-1}}\big|\alpha(\alpha-1)\psi_t^2+\alpha\psi\psi_{tt}-\alpha(\alpha-1)|\nabla\psi|^2-\alpha\psi\Delta\psi\big|^\frac{q}{q-1}\,d\mu dt\right)^{\frac{q-1}q}\\ &
\times \left(\int_{S_K} |u|^q\psi^\alpha V\,d\mu dt \right)^{\frac 1 q}\,.
\end{aligned}
  \end{equation}
\end{lemma}
\begin{proof}
We use Definition \ref{defsol}, with the test function $\varphi=\psi^{\alpha}$ and $\psi\equiv 1$ in $K$. By the same arguments as in the proof of Lemma \ref{lemma1}, we get
\[\int_0^{\infty}\int_M |u|^q\psi^{\alpha} V\, d\mu dt \leq  \int_{S_K} |u||(\psi^{\alpha})_{tt}-\Delta(\psi^{\alpha})| d\mu dt
\]
\[+\int_M [u_0(x) (\psi^{\alpha})_t(x,0)-u_1(x)\psi^{\alpha}(x,0)]\,d\mu
\]
\[\leq \left(\int_{S_K} |u|^q \psi^{\alpha} V\, d\mu dt\right)^{\frac 1 q}\left(\int_{S_K}|(\psi^{\alpha})_{tt} - \Delta(\psi^{\alpha})|^{\frac{q}{q-1}}(V \psi^{\alpha})^{-\frac{1}{q-1}}\, d\mu dt\right)^{\frac {q-1}q}\]
\[ + \int_M [u_0(x) (\psi^{\alpha})_t(x,0)- u_1(x)\psi^{\alpha}(x,0)]\,d\mu
\leq - \int_M u_1(x,0)\psi^\alpha(x,0)-u_0(x,0)(\psi^\alpha(x,0))_t\,d\mu \]
\[+\left(\int_{S_K} V^{-\frac{1}{q-1}}\psi^{\alpha-\frac{2q}{q-1}}\big|\alpha(\alpha-1)\psi_t^2+\alpha\psi\psi_{tt}-\alpha(\alpha-1)|\nabla\psi|^2-\alpha\psi\Delta\psi\big|^\frac{q}{q-1}\,d\mu dt\right)^{\frac{q-1}q}\]
\[
\times \left(\int_{S_K} |u|^q\psi^\alpha V\,d\mu dt \right)^{\frac 1 q}\,,
\]
which is the thesis.
\end{proof}

\section{Proof of the main results}\label{proof}
\begin{proof}[Proof of Theorem \ref{teo1}] Let us begin by assuming that $\sigma\in [-2, 2]$. Consider the family of test functions $\{\phi_R\}_{R\geq 1}$ provided by Proposition \ref{propSB}. Next, let
$\eta_R\in C^{2}([0, \infty))$ with the following properties
\begin{itemize}
\item[$(i)$] $0\leq \eta_R\leq 1$, $\eta_R\equiv 1$ in $[0, R]$;
\item[$(ii)$] $\operatorname{supp} \eta_R\subset [0, \gamma R)$ with the same $\gamma$ as in Proposition \ref{propSB};
\item[$(iii)$] $-\frac{C_1}{R}\eta_R'\leq 0$ for some $C_1>0$ independent of $R$;
 \item[$(iv)$] $\left|\eta_R''\right|\leq \frac{C_1}{R^2}$ for some $C_1>0$ independent of $R$.
\end{itemize}
Fix any $\alpha>\frac{2q}{q-1}$. Now, for every $R>1$ we apply Lemma \ref{lemma1} with
\begin{equation}\label{eq13}
\psi_R(x,t)\equiv \psi(x,t):= \phi_R(x)\eta_R(t), \qquad x\in M,\, t\geq 0\,.
\end{equation}
Note that, in view of the properties of the functions $\phi_R$ and $\eta_R$, for some constant $\tilde C>0$, there holds
\begin{equation}\label{eq7}
\begin{aligned}
\left|\alpha(\alpha-1)(\psi_t)^2+\alpha\psi\psi_{tt}-\alpha(\alpha-1)|\nabla\psi|^2-\alpha\psi\Delta\psi \right|^{\frac{q}{q-1}}\\
\leq \tilde C\left(\frac{1}{R^2}+\frac{1}{R^{1+\frac{\sigma}{2}}}\right)^{\frac{q}{q-1}}\leq \tilde C R^{-\left(1+\frac{\sigma}{2}\right)\frac{q}{q-1}}
\end{aligned}
\end{equation}
for every $R>1$.
It is direct to see that for all $x\in M$
\begin{equation}\label{eq9}
(\psi^{\alpha})_t(x, 0)=0\,.
\end{equation}
By \eqref{eq3}, \eqref{eq7}, \eqref{eq9},
\begin{equation}\label{eq8}
\begin{aligned}
\int_0^R\int_{B_R(o)} |u|^q V\, d\mu dt \leq  \int_0^{\infty}\int_{M} |u|^q \psi^{\alpha}V\, d\mu dt
\\ \leq \tilde C R^{-\left(1+\frac{\sigma}{2}\right)\frac{q}{q-1}}\int_0^{\gamma R}\int_{B_{\gamma R}(o)} V^{-\frac{1}{q-1}}\, d\mu dt \\
- \frac{q}{q-1}\int_{B_{\gamma R}(o)} u_1\phi_R^{\alpha}(x)\, d\mu\,.
\end{aligned}
\end{equation}
Now, by condition \eqref{eq11}, we can infer that, for all $R>\frac{R_0}{\gamma},$
\begin{equation}\label{eq16}
  \int_0^{\gamma R}\int_{B_{\gamma R}(o)} V^{-\frac{1}{q-1}}\, d\mu dt \leq \bar C \gamma^{\left(1+\frac{\sigma}{2}\right)\frac{q}{q-1}} R^{\left(1+\frac{\sigma}{2}\right)\frac{q}{q-1}}\,.
\end{equation}
Moreover,
\begin{equation}\label{eq15}
- \int_{B_{\gamma R}(o)} u_1\phi_R^{\alpha}(x)\, d\mu \leq -\int_{B_R(o)}u_1^+ d\mu + \int_{B_{\gamma R}(o)} u_1^-\, d\mu \mathop{\longrightarrow}_{R\to+\infty} - \int_M u_1\, d\mu \leq 0\,.
\end{equation}
Passing to the limit as $R\to\infty$ in \eqref{eq8} we obtain
\begin{equation}\label{eq12}
\int_0^{\infty}\int_{M} |u|^q V\, d\mu dt \leq C
\end{equation}
for some constant $C\in (0, \infty)$.
For every $R>R_0$, we can use Lemma \ref{lemma2} with $\psi\equiv \psi_R$ given by \eqref{eq13} and $K\equiv K_R=[0, R]\times\overline{B_R(o)}$.
This yields
\begin{equation}\label{eq14}
\begin{aligned}
&\int_{K} |u|^q V\,d\mu dt\leq - \int_Mu_1(x,0)\psi^\alpha(x,0)-u_0(x,0)(\psi^\alpha(x,0))_t\,d\mu \\ &+\left(\int_{S_K} V^{-\frac{1}{q-1}}\psi^{\alpha-\frac{2q}{q-1}}\big|\alpha(\alpha-1)\psi_t^2+\alpha\psi\psi_{tt}-\alpha(\alpha-1)|\nabla\psi|^2-\alpha\psi\Delta\psi\big|^\frac{q}{q-1}\,d\mu dt\right)^{\frac{q-1}q}\\ &
\times \left(\int_{S_K} |u|^q\psi^\alpha V\,d\mu dt \right)^{\frac 1 q}\,.
\end{aligned}
  \end{equation}
In view of \eqref{eq7}, \eqref{eq9}, \eqref{eq8} and \eqref{eq16} we obtain
\begin{equation}\label{eq17}
\int_{K} |u|^q V\,d\mu dt\leq - \int_Mu_1(x,0)\psi^\alpha(x,0)\,d\mu + C \left(\int_{S_K} |u|^qV\,d\mu dt \right)^{\frac 1 q}\,,
\end{equation}
for some $C>0$ independent of $R$. Letting $R\to\infty$ in \eqref{eq17}, by \eqref{eq15} and \eqref{eq12} we have
\[\int_0^{\infty}\int_M |u|^q V\, d\mu dt \leq 0\,. \]
Since $V>0$ a.e. in $M\times(0, \infty)$, we deduce that $u\equiv 0$ a.e. in $M\times(0, \infty)$. This completes the proof in the case $\sigma\in [-2, 2]$.

\smallskip

Now, note that if $\sigma>2$, then \eqref{eq4} is also satisfied for $\sigma =2$. Then the conclusion follows from the case $\sigma=2$ already treated above.
\end{proof}

\begin{proof}[Proof of Theorem \ref{teo2}] Let $\zeta\in C^2([0, \infty))$ with $\zeta\geq 0,\,\zeta'\leq 0, \, \operatorname{supp} \zeta\subseteq [0, 2],\, \zeta\equiv 1$ in $[0, 1]$.
In view of \eqref{eq34}, \eqref{eq41} holds with $\phi\in \mathcal A$, $\phi(r)=\exp\{B r^{1-\frac{\sigma}{2}}\}$ for every $r>1$, for some $B>0$. So, due to \eqref{eq42},
\[0\leq \mathcal{F}(r, \theta)\leq  C (1+r)^{-\frac{\sigma}{2}}\quad \textrm{for all}\;\; r>\bar R\,,\]
for some $C>0, \bar R>0$. Since $M$ is a Cartan-Hadamard manifold, $\operatorname{Cut}(o)=\emptyset$. Hence, $x\mapsto r(x)$ is a function of class $C^{2}$ in $M\setminus\{o\}$. Thus, \eqref{eq38} holds and it yields
\begin{equation}\label{eq43}
\begin{aligned}
\left|\Delta\left[\zeta\left(\frac{r(x)}{R}\right)\right]\right| & \leq \frac{1}{R^2}\left|\zeta''\left(\frac{r(x)}{R}\right)\right| +
C \frac{n-1}{R}(1+r)^{-\frac{\sigma}{2}}\left|\zeta'\left(\frac{r(x)}{R}\right)\right|\\
& \leq  C (1+R)^{-\frac{\sigma}{2}-1} \quad \textrm{for all} \;\; x\in M\,.
\end{aligned}
\end{equation}
For every $R>1$ let
\[\psi(x,t)\equiv \psi_R(x, t):=\zeta\left(\frac{r(x)}{R}\right)\zeta\left(\frac{t}{R}\right)\quad x\in M, t>0\,.\]
Note that, in view of the properties of \eqref{eq43} and the properties of the function $\zeta$, for some constant $\tilde C>0$, there holds
\begin{equation}\label{eq44}
\begin{aligned}
\left|\alpha(\alpha-1)(\psi_t)^2+\alpha\psi\psi_{tt}-\alpha(\alpha-1)|\nabla\psi|^2-\alpha\psi\Delta\psi \right|^{\frac{q}{q-1}}\\
\leq \tilde C\left(\frac{1}{R^2}+\frac{1}{R^{1+\frac{\sigma}{2}}}\right)^{\frac{q}{q-1}}\leq \tilde C R^{-\left(1+\frac{\sigma}{2}\right)\frac{q}{q-1}}
\end{aligned}
\end{equation}
for every $R>1$.
By arguing as in the proof of Theorem \ref{teo1} with \eqref{eq7} replaced by \eqref{eq44} the conclusion follows.
\end{proof}


\begin{thebibliography}{999}


\bibitem{AMR} L. J. Al\'ias, P. Mastrolia,  M. Rigoli, {\it Maximum principles and geometric applications}, Springer, 2016.

\bibitem{BPT} C. Bandle, M.A. Pozio, A. Tesei, {\it The Fujita Exponent for the Cauchy Problem in the
Hyperbolic Space}, J. Diff. Eq. {\bf 251} (2011), 2143--2163 \,.

\bibitem{BS} D. Bianchi, A. Setti, {\it Laplacian cut-offs, porous and fast diffusion on manifolds and other applications}, Calc. Var. Part. Diff. Eq. {\bf 57}:4 (2018).

\bibitem{survey}
K. Deng, H.A.\ Levine, {\it The role of critical exponents in blow-up theorems: the sequel},
J. Math. Anal. Appl. {\bf 243} (2000), 85--126.


\bibitem{fuji}
H. Fujita, {\it On the blowing up of solutions of the Cauchy problem $u_t=\Delta+u^{1+\alpha}$},
J. Fac. Sci. Univ. Tokio Sec {\bf 13} (1966), 109--124.

\bibitem{Grig}
A. Grigoryan, \emph{ Analytic and geometric background of
recurrence and non-explosion of the Brownian motion on Riemannian
manifolds}, \newblock Bull. Amer. Math. Soc. {\bf 36} (1999),
135--249.


\bibitem{GrigKond} A. Grigor'yan, V. A. Kondratiev, {\it On the existence of positive solutions of semilinear elliptic inequalities on Riemannian manifolds}, In Around the research of Vladimir Maz'ya. II, volume 12 of Int. Math. Ser. (N. Y.), pages
203--218. Springer, New York, 2010.

\bibitem{GrigS} A. Grigor'yan, Y. Sun, {\it On non-negative solutions of the inequality $\Delta u + u^\sigma \leq 0$ on Riemannian manifolds}\,, Comm. Pure Appl. Math.
{\bf 67} (2014), 1336--1352\,.

\bibitem{haya}
K. Hayakawa, {\it On the nonexistence of global solutions of some semilinear
	parabolic equations},
Proc. Japan Acad Sci. {\bf 49} (1973), 503--505.

\bibitem{kato}
T. Kato, {\it Blow-up of solutions of some nonlinear hyperbolic equations},
Comm. Pure Appl. Math {\bf 33} (1980), 501--505

\bibitem{MMP1} P. Mastrolia, D.D. Monticelli, F. Punzo, emph{Nonexistence results for elliptic differential inequalities with a potential on Riemannian manifolds}, Calc. Var.Part. Diff. Eq. {\bf 54} (2015), 1345-1372;

\bibitem{MMP2}
P. Mastrolia, D. Monticelli, F.Punzo,
{\it Nonexistence of solutions to parabolic differential inequalities with a potential on Riemannian manifolds}, Math Ann. {\bf 367} (2017), 929--963.

\bibitem{MP}
E. Mitidieri, S. Pohozaev, {\it Nonexistence of weak solutions for some degenrate and singular hyperbolic problems on $\mathbb R^n$}, Proc. Steklov Inst. Math. {\bf 252} (2001), 1--19\,.

\bibitem{Pu1} F. Punzo, {\it Blow-up of solutions to semilinear parabolic equations on Riemannian manifolds with negative sectional curvature}, J. Math. Anal. Appl. {\bf 387} (2012) 815--827;


\bibitem{shae}
J. Shaeffer, {\it The equation $u_{tt}-\Delta u=|u|^p$ for the critical value of $p$}, Proc. Royal. Soc Edinburgh, {\bf 101} (1985), 31--44.



\bibitem{Ru} Q. Ru, {\it A nonexistence result for a nonlinear wave equation with damping on a Riemannian manifold}, Boundary Value Problems {\bf 198} (2016)\,.

\bibitem{zhang}
Q.S.\ Zhang,
{\it Blow-up results for nonlinear parabolic equations on manifolds},
Duke Math. J.\ {\bf 97} (1999), 515--539.

\end{thebibliography}
\end{document}